\documentclass[11pt,reqno]{amsart}
\usepackage[a4paper,margin=27mm]{geometry}
\usepackage[T1]{fontenc}
\usepackage{lmodern}
\usepackage{amsmath,amssymb,amsthm}
\usepackage{microtype}
\usepackage{cite}
\usepackage{indentfirst}
\usepackage[hidelinks]{hyperref}
\hypersetup{
  pdftitle={A proof of the Letzter\textendash Morrison conjecture on reduced multisets},
  pdfsubject={The exact value of the reduced-multiset parameter in dimension two},
  pdfkeywords={subsequence sums, reduced multisets, additive bases, additive combinatorics}
}
\newtheorem{theorem}{Theorem}[section]
\newtheorem{conjecture}[theorem]{Conjecture}
\newtheorem{lemma}[theorem]{Lemma}
\newtheorem{corollary}[theorem]{Corollary}
\theoremstyle{remark}
\newtheorem{remark}[theorem]{Remark}
\numberwithin{equation}{section}
\title{A proof of the Letzter--Morrison conjecture on reduced multisets}
\author{Wenkai Yang \and Chuhan Lei}

\address{Center for Combinatorics \\ Nankai University, Tianjin 300071, P.R. China}
\email{9820250116@nankai.edu.cn}

\address{School of Mathematical Sciences \\ Tiangong University \\ Tianjin 300387, P.R. China}
\email{863100953@qq.com}

\date{}

\begin{document}

\begin{abstract}
A finite multiset over an abelian group is reduced if every proper submultiset has fewer distinct subsequence sums.
Letzter and Morrison introduced reduced multisets in 2023 to study zero-sum cycles in complete digraphs.
For a prime $p$, let $h_p(2)$ be the maximum size of a reduced multiset over $C_p^2$.
Letzter and Morrison~\cite{LM} established the lower bound $h_p(2)\ge2p-2$ and conjectured that $h_p(2)\le2p+c$ for some absolute constant $c$.
In this paper, we use two different methods to prove that $h_p(2)=2p-2$ for every prime $p$, and formed this construction.
\end{abstract}
\maketitle

\section{Introduction}\label{sec:introduction}

Let $G$ be a finite abelian group.
We regard a finite multiset over $G$ as an unordered sequence.
For such a sequence $S$, let $\Sigma_0(S)$ denote the set of all subsequence sums, including the empty sum.
Following Letzter and Morrison~\cite{LM}, we call $S$ \emph{reduced} if
\begin{equation}\label{eq:reduced}
 \Sigma_0(Sa^{-1})\subsetneq\Sigma_0(S)\qquad\text{for every }a\mid S,
\end{equation}
where $a\mid S$ means that $S=aT$ for some sequence $T$, in which case $Sa^{-1}=T$.
For a prime $p$, let $C_p$ denote the cyclic group of order $p$.
For each integer $k\ge1$, let $h_p(k)$ denote the maximum length of a reduced sequence over $C_p^k$.

The additive invariant $h_p(k)$ also arises in the study of zero-sum cycles in complete digraphs.
Alon and Krivelevich~\cite{AK} studied how many vertices a complete digraph with arc weights in a cyclic group must have to guarantee a directed cycle of total weight zero.
M\'esz\'aros and Steiner~\cite{MS} obtained bounds for finite abelian groups.
For weights in $C_p^k$, Letzter and Morrison~\cite[Lemma~4.1]{LM} used $h_p(10k)$ to bound the required number of vertices.
Their bound on $h_p(k)$ is related to the problem of covering $C_p^k$ by subsequence sums of a multiset union of bases, studied by Jaeger, Linial, Payan and Tarsi~\cite{JLPT} and by Alon, Linial and Meshulam~\cite{ALM}.
Arsovski~\cite{Arsovski} studied related questions for generating sets of finite abelian groups.

Letzter and Morrison~\cite[Theorem~2.1]{LM} proved that $h_p(k)\le k(p-1)(1+\log_2 k)$.
For $k=2$ and $p\ge7$, they proved the sharper bound~\cite[Lemma~5.2]{LM}
\[
 h_p(2)<\frac{5}{2}(p-1).
\]

Furthermore, they proposed the following conjecture.

\begin{conjecture}\label{conj:LM}\textup{\cite[Conjecture~5.3]{LM}}
There exists an absolute constant $c$ such that for every prime $p$,
\[
 h_p(2)\le2p+c.
\]
\end{conjecture}

In this paper, we prove the following theorem.

\begin{theorem}\label{thm:main}
For every prime $p$,
\[
 h_p(2)=2p-2.
\]
\end{theorem}

Theorem~\ref{thm:main} settles the conjecture of Letzter and Morrison with the best possible constant $c=-2$.
The invariant $h_p(2)$ is closely related to the additive basis problem for $C_p^2$, which has been studied extensively~\cite{MO,MW,PengI}.
For every regular sequence $S=g_1\cdots g_{2p-1}$ over $C_p^2$, Theorem~\ref{thm:deletion} shows that $\Sigma_0(Sg_i^{-1})=C_p^2$ for at least $p$ indices $i$.
By Corollary~\ref{cor:preserve}, if $S$ is a finite sequence over $C_p^2$ and $Z\mid S$ satisfies $0\nmid Z$ and $|Z|\le p-1$, then there exists a subsequence $R$ with $Z\mid R\mid S$, $|R|\le2p-2$ and $\Sigma_0(R)=\Sigma_0(S)$.

In Section~\ref{sec:preliminaries} we establish basic properties of reduced sequences and prove that every reduced sequence is regular.
Section~\ref{sec:proof} contains the proofs of Theorems~\ref{thm:main} and~\ref{thm:deletion} and Corollary~\ref{cor:preserve}.
Section~\ref{sec:remarks} contains an alternative proof of the upper bound.

\section{Preliminaries}\label{sec:preliminaries}

All groups are finite abelian groups with additive notation.
A sequence over $G$ has the form
\[
 S=g_1\cdot\ldots\cdot g_\ell
   =\prod_{g\in G}g^{\mathsf v_g(S)},
\]
where $\mathsf v_g(S)$ is the multiplicity of $g$ in $S$ and $|S|=\ell$ is the length of $S$.
The product $SR$ denotes concatenation, and $1$ denotes the empty sequence.
We identify each $a\in G$ with the sequence of length one consisting of $a$.

A sequence $U$ is a subsequence of $S$, denoted by $U\mid S$, if $\mathsf v_g(U)\le\mathsf v_g(S)$ for every $g\in G$.
If $U\mid S$, then $SU^{-1}$ denotes the unique sequence such that $S=U(SU^{-1})$.
In particular, $\mathsf v_g(SU^{-1})=\mathsf v_g(S)-\mathsf v_g(U)$ for every $g\in G$.

For $S=g_1\cdot\ldots\cdot g_\ell$, let
\[
 \sigma(S)=g_1+\cdots+g_\ell,
 \qquad \Sigma_0(S)=\{\sigma(U):U\mid S\}.
\]
We use the convention $\sigma(1)=0$, so $\Sigma_0(1)=\{0\}$.
Let $\Sigma(S)=\{\sigma(T):1\ne T\mid S\}$ be the set of nonempty subsequence sums.
Thus $\Sigma_0(S)=\Sigma(S)\cup\{0\}$ and $\Sigma(1)=\emptyset$.

For subsets $B,D\subseteq G$ and $x\in G$, let $B+D=\{b+d:b\in B,\ d\in D\}$, $B+x=\{b+x:b\in B\}$, and $x-B=\{x-b:b\in B\}$.
For sequences $S,R$ over $G$ and $x\in G$, we have
\begin{equation}\label{eq:union}
 \begin{aligned}
 \Sigma_0(SR)&=\Sigma_0(S)+\Sigma_0(R),\\
 \Sigma_0(Sx)&=\Sigma_0(S)\cup(\Sigma_0(S)+x).
 \end{aligned}
\end{equation}

For a subgroup $H\le G$, let $S_H=\prod_{g\in H}g^{\mathsf v_g(S)}$.
The sequence $S$ is called \emph{regular} if $|S_H|\le |H|-1$ for every proper subgroup $H<G$~\cite{PQL}.
A sequence $S$ is an \emph{additive basis} of $G$ if $\Sigma(S)=G$.
Let $\mathsf c_0(G)$ denote the smallest integer $\ell$ such that every regular sequence $S$ over $G$ with $|S|\ge\ell$ is an additive basis~\cite{GHQQZ,QLWZ}.
In particular, an additive basis must have a nonempty zero-sum subsequence.

Lemma~\ref{lem:elementary}(1) is due to Letzter and Morrison~\cite[Observation~2.2(i), (iii)]{LM}.

\begin{lemma}\label{lem:elementary}
Let $G$ be a finite abelian group and let $S$ be a sequence over $G$.
\begin{enumerate}
\item If $S$ is reduced and $T\mid S$, then $0\nmid S$ and $T$ is reduced.
\item If $S$ is reduced, then $|\Sigma_0(S)|\ge |S|+1$.
\item We have
\begin{equation}\label{eq:reflection}
 \Sigma_0(S)=\sigma(S)-\Sigma_0(S).
\end{equation}
Moreover, $\sigma(S)-(G\setminus\Sigma_0(S))=G\setminus\Sigma_0(S)$.
\end{enumerate}
\end{lemma}

\begin{proof}
\textup{(1)} If $S=0R$ for some sequence $R$, then $\Sigma_0(S)=\Sigma_0(R)$, contrary to~\eqref{eq:reduced} with $a=0$.
Thus $0\nmid S$.

Assume to the contrary that $T$ is not reduced.
Then there exists $a\mid T$ such that $\Sigma_0(Ta^{-1})=\Sigma_0(T)$.
Let $R=ST^{-1}$, so that $S=TR$ and $Sa^{-1}=(Ta^{-1})R$.
By~\eqref{eq:union},
\[
 \Sigma_0(Sa^{-1})=\Sigma_0(Ta^{-1})+\Sigma_0(R)=\Sigma_0(T)+\Sigma_0(R)=\Sigma_0(S).
\]
Since $a\mid S$ and $S$ is reduced, we have $\Sigma_0(Sa^{-1})\subsetneq\Sigma_0(S)$, a contradiction.

\textup{(2)} If $S=1$, then $|\Sigma_0(S)|=1=|S|+1$.
Otherwise, let $S=g_1\cdot\ldots\cdot g_\ell$, where $\ell=|S|\ge1$.
For $0\le i\le\ell$, let $S_i=\prod_{j=1}^i g_j$, with $S_0=1$.
By \textup{(1)}, each $S_i$ is reduced since $S_i\mid S$.
For $1\le i\le\ell$, we have $S_{i-1}=S_i g_i^{-1}$ and hence $\Sigma_0(S_{i-1})\subsetneq\Sigma_0(S_i)$.
Since $|\Sigma_0(S_0)|=1$, it follows that $|\Sigma_0(S)|\ge1+\ell=|S|+1$.

\textup{(3)} For every $T\mid S$, the complementary subsequence $ST^{-1}$ satisfies $\sigma(ST^{-1})=\sigma(S)-\sigma(T)$.
The map $T\mapsto ST^{-1}$ is its own inverse on the subsequences of $S$, so~\eqref{eq:reflection} follows.
By~\eqref{eq:reflection}, the bijection $x\mapsto\sigma(S)-x$ of $G$ maps $\Sigma_0(S)$ onto itself and hence also maps $G\setminus\Sigma_0(S)$ onto itself.
\end{proof}

\begin{lemma}\label{lem:line-bound}
Let $G$ be a finite abelian group and let $S$ be a reduced sequence over $G$.
Then $S$ is regular.
More precisely,
\[
 |S_H|\le |H|-1\qquad\text{for every subgroup $H\le G$}.
\]
In particular, if $G=C_p^2$ for a prime $p$, then $|S_H|\le p-1$ for every subgroup $H$ of order $p$.
\end{lemma}

\begin{proof}
For every subgroup $H\le G$, Lemma~\ref{lem:elementary}(1) shows that $S_H$ is reduced since $S_H\mid S$.
By Lemma~\ref{lem:elementary}(2) and $\Sigma_0(S_H)\subseteq H$,
\[
 |S_H|+1\le |\Sigma_0(S_H)|\le |H|.
\]
Thus $|S_H|\le |H|-1$ for all $H\le G$, so $S$ is regular.
If $|H|=p$, then $|S_H|\le p-1$.
\end{proof}

Peng~\cite{PengI} proved that $\mathsf c_0(C_p^2)=2p-1$ for every prime $p$ (see also~\cite[p.~266]{GHQQZ}).
The proof of Theorem~\ref{thm:main} uses Peng's covering theorem~\cite[Theorem~2]{PengI} in the following form, stated by Matom\"aki~\cite[Lemma~14]{M}.

\begin{theorem}\label{thm:peng}\textup{\cite{PengI}}
Let $p$ be a prime, let $G=C_p^2$, and let $S$ be a sequence over $G$ with $|S|=2p-1$ and $0\nmid S$.
If $|S_H|\le p$ for every subgroup $H$ of order $p$, then $\Sigma_0(S)=G$.
\end{theorem}

\section{Proof of the main theorem}\label{sec:proof}

Throughout this section, $p$ is a prime and $G=C_p^2$.
For odd $p$ and $x\in G$, let $x/2=\frac{p+1}{2}x$.
Since $px=0$, we have $2\bigl(\frac{p+1}{2}x\bigr)=(p+1)x=x$.
Thus multiplication by $(p+1)/2$ is the inverse of multiplication by $2$, and $x/2$ is the unique $y\in G$ such that $2y=x$.

\begin{theorem}\label{thm:deletion}
Let $p$ be a prime and let $S=g_1\cdots g_{2p-1}$ be a regular sequence over $G=C_p^2$.
If $p=2$, then $\Sigma_0(Sg_i^{-1})=G$ for every $i$.
If $p$ is odd, then either $\Sigma_0(Sa^{-1})=G$ for every $a\mid S$, or there exists $u\mid S$ such that
\[
 \Sigma_0(Sa^{-1})\ne G\quad\Longleftrightarrow\quad a\in\{u,-u\}
 \qquad\text{for every }a\mid S.
\]
In particular, at least $p$ indices $i$ satisfy $\Sigma_0(Sg_i^{-1})=G$, and this bound is sharp for every odd prime $p$.
\end{theorem}

\begin{proof}
Since $S$ is regular, $0\nmid S$ and, for every subgroup $H$ of order $p$,
\begin{equation}\label{eq:line-bound}
 |S_H|\le p-1.
\end{equation}
If $p=2$, let $e_1,e_2$ be a basis of $G$.
By~\eqref{eq:line-bound}, $S=e_1e_2(e_1+e_2)$ is squarefree.
For every $a\mid S$, the two elements of $Sa^{-1}$ form a basis of $G$, so $\Sigma_0(Sa^{-1})=G$.
We now assume that $p$ is odd.
By~\eqref{eq:line-bound} and Theorem~\ref{thm:peng}, $\Sigma_0(S)=G$.

We first show that $|\Sigma_0(S)\setminus\Sigma_0(Sa^{-1})|\le1$ for every $a\mid S$, following the argument of Matom\"aki~\cite[Lemma~15]{M}.
Assume to the contrary that there exists $a\mid S$ such that $|\Sigma_0(S)\setminus\Sigma_0(Sa^{-1})|\ge2$.
Then there exist distinct $u,v\in\Sigma_0(S)\setminus\Sigma_0(Sa^{-1})$.
Since $0\nmid S$ and $u-v\ne0$, we have $0\nmid Sa^{-1}\cdot(u-v)$ and $|Sa^{-1}\cdot(u-v)|=2p-1$.
For every subgroup $H$ of order $p$, we have $|(Sa^{-1})_H|\le |S_H|$ and $|(u-v)_H|\le1$, so by~\eqref{eq:line-bound},
\[
 |(Sa^{-1}\cdot(u-v))_H|=|(Sa^{-1})_H|+|(u-v)_H|\le |S_H|+1\le p.
\]
By Theorem~\ref{thm:peng}, $\Sigma_0(Sa^{-1}\cdot(u-v))=G$.
Since $u\notin\Sigma_0(Sa^{-1})$ and $u-(u-v)=v\notin\Sigma_0(Sa^{-1})$, equation~\eqref{eq:union} implies
\[
 u\notin\Sigma_0(Sa^{-1})\cup\bigl(\Sigma_0(Sa^{-1})+(u-v)\bigr)=\Sigma_0(Sa^{-1}\cdot(u-v)),
\]
a contradiction.

If $a\mid S$ and $\Sigma_0(Sa^{-1})\ne G$, let $m_a$ be the unique element of $G\setminus\Sigma_0(Sa^{-1})$.
By Lemma~\ref{lem:elementary}(3) for $Sa^{-1}$, the map $x\mapsto\sigma(Sa^{-1})-x$ maps $\{m_a\}$ onto itself.
Hence $\sigma(S)-a-m_a=m_a$, and therefore
\begin{equation}\label{eq:single-deletion}
 m_a=\frac{\sigma(S)-a}{2}\qquad\text{whenever }a\mid S\text{ and }\Sigma_0(Sa^{-1})\ne G.
\end{equation}

We next show that if $a,b\mid S$ and both $\Sigma_0(Sa^{-1})$ and $\Sigma_0(Sb^{-1})$ differ from $G$, then $b\in\{a,-a\}$.
Assume to the contrary that there exist such $a,b$ with $b\notin\{a,-a\}$.
Since $a\ne b$, we have $ab\mid S$.
Let $T=Sa^{-1}b^{-1}$ and $c=\sigma(T)/2=(\sigma(S)-a-b)/2$.
By~\eqref{eq:single-deletion}, \eqref{eq:union} and $Sa^{-1}=Tb$, we have
\[
 m_a=c+b/2\notin\Sigma_0(T)\cup(\Sigma_0(T)+b).
\]
Thus $c+b/2=m_a\notin\Sigma_0(T)$ and $c-b/2=m_a-b\notin\Sigma_0(T)$, since $m_a\notin\Sigma_0(T)+b$.
Similarly, using $Sb^{-1}=Ta$, we obtain $c+a/2=m_b\notin\Sigma_0(T)$ and $c-a/2=m_b-a\notin\Sigma_0(T)$.
Since $\Sigma_0(S)=G$, it follows that
\begin{equation}\label{eq:four-points}
 \{c+a/2,c-a/2,c+b/2,c-b/2\}\subseteq\Sigma_0(S)\setminus\Sigma_0(T).
\end{equation}
These four elements are distinct because $p$ is odd, $a,b\ne0$, and $a\notin\{b,-b\}$.

For $x=(a+b)/2$ and $y=(a-b)/2$, we show that $\Sigma_0(Txy)=G$.
Since $b\notin\{a,-a\}$, we have $x,y\ne0$, so $0\nmid Txy$ and $|Txy|=|T|+2=2p-1$.
If $\langle a\rangle=\langle b\rangle$, then $a,b,x,y$ all belong to this subgroup and are nonzero.
Thus $|(Txy)_H|=|S_H|$ for every subgroup $H$ of order $p$, since distinct subgroups of order $p$ intersect only in $\{0\}$.
If $\langle a\rangle\ne\langle b\rangle$, no subgroup $H$ of order $p$ contains both $x$ and $y$, since it would also contain $a=x+y$ and $b=x-y$.
In this case $|(Txy)_H|\le|T_H|+1\le|S_H|+1\le p$ by~\eqref{eq:line-bound}.
Thus in both cases $\Sigma_0(Txy)=G$ by Theorem~\ref{thm:peng}.

We show that $m_b\notin\Sigma_0(Txy)$.
Since $m_b=c+a/2$, the definitions of $x$ and $y$ imply
\[
 \begin{aligned}
 m_b&=c+a/2,       & m_b-x&=c-b/2,\\
 m_b-y&=c+b/2,     & m_b-x-y&=c-a/2.
 \end{aligned}
\]
By~\eqref{eq:four-points}, $\{m_b,m_b-x,m_b-y,m_b-x-y\}\cap\Sigma_0(T)=\emptyset$.
By~\eqref{eq:union}, $\Sigma_0(Txy)=\Sigma_0(T)+\{0,x,y,x+y\}$.
Thus $m_b\notin\Sigma_0(Txy)$, although $\Sigma_0(Txy)=G$, a contradiction.

If $\Sigma_0(Sa^{-1})=G$ for every $a\mid S$, then all $2p-1$ indices satisfy $\Sigma_0(Sg_i^{-1})=G$.
Otherwise there exists $u\mid S$ such that $\Sigma_0(Su^{-1})\ne G$, and every $a\mid S$ with $\Sigma_0(Sa^{-1})\ne G$ belongs to $\{u,-u\}$.
If $(-u)\mid S$, let $V=S u^{-1}(-u)^{-1}$.
By~\eqref{eq:union}, $\Sigma_0(S(-u)^{-1})=\Sigma_0(Vu)=u+\Sigma_0(V(-u))=u+\Sigma_0(Su^{-1})\ne G$.
Thus $\Sigma_0(Sa^{-1})\ne G$ if and only if $a\in\{u,-u\}$ for $a\mid S$.
By~\eqref{eq:line-bound}, $\mathsf v_u(S)+\mathsf v_{-u}(S)\le |S_{\langle u\rangle}|\le p-1$, so at least $(2p-1)-(p-1)=p$ indices satisfy $\Sigma_0(Sg_i^{-1})=G$.

For sharpness, let $e_1,e_2$ be a basis of $G$ and let $H=\langle e_1\rangle$.
Following Martin, Peilloux and Wong~\cite[Examples~4.2 and~4.4]{MPW} in dimension two, let
\begin{equation}\label{eq:deletion-example}
 Z=\prod_{j=0}^{p-1}(je_1+e_2),\qquad W=e_1^{\,p-1}Z.
\end{equation}
For $0\le i,j<p$, the inclusion $ie_1+e_2\in\langle je_1+e_2\rangle$ implies $i=j$ by comparison of the two coordinates.
Thus the $p$ elements of $Z$ generate distinct subgroups of order $p$, all different from $H$, so $W$ is regular and $|W|=2p-1$.
For every $U\mid Z$, the second coordinate of $\sigma(U)$ is $|U|$ modulo $p$.
Thus $\sigma(U)\in H$ only for $U=1$ or $U=Z$, and $\sigma(Z)=\frac{p(p-1)}2e_1+pe_2=0$ since $p$ is odd.
It follows that $\Sigma_0(We_1^{-1})\cap H=\{je_1:0\le j\le p-2\}$, so $-e_1\notin\Sigma_0(We_1^{-1})$.
For every $b\mid Z$, the sums of subsequences of $Zb^{-1}$ of lengths $0,\ldots,p-1$ represent all cosets of $H$.
Since $\Sigma_0(e_1^{\,p-1})=H$, we have $\Sigma_0(Wb^{-1})=H+\Sigma_0(Zb^{-1})=G$.
Thus the $p$ elements of $Z$ are precisely the $b\mid W$ for which $\Sigma_0(Wb^{-1})=G$.
\end{proof}

\begin{proof}[Proof of Theorem~\ref{thm:main}]
	
\textit{Lower bound.} Let $e_1,e_2$ be a basis of $G$, and let $S_0=e_1^{\,p-1}e_2^{\,p-1}$.
Since $G=\langle e_1\rangle\oplus\langle e_2\rangle$, we have
\[
 \Sigma_0(S_0)=\{ie_1+je_2:0\le i,j\le p-1\}=G,
\]
and the two coefficients are unique modulo $p$.
For $S_0e_1^{-1}=e_1^{p-2}e_2^{p-1}$, the first coefficient ranges from $0$ to $p-2$, so $-e_1\in\Sigma_0(S_0)\setminus\Sigma_0(S_0e_1^{-1})$.
Similarly, $-e_2\in\Sigma_0(S_0)\setminus\Sigma_0(S_0e_2^{-1})$.
Thus $S_0$ is reduced and $h_p(2)\ge |S_0|=2p-2$.

\textit{Upper bound.} Assume to the contrary that $S$ is a reduced sequence over $G$ of length $2p-1$.
By Lemma~\ref{lem:line-bound}, $S$ is regular.
By Theorem~\ref{thm:deletion}, there exists $a\mid S$ such that $\Sigma_0(Sa^{-1})=G=\Sigma_0(S)$, contrary to~\eqref{eq:reduced}.
Hence no reduced sequence over $G$ has length $2p-1$.

It follows that $h_p(2)\le2p-2$.
Indeed, for any sequence $S$ over $G$ with $|S|\ge2p-1$, there exists $T\mid S$ with $|T|=2p-1$.
Since no sequence of length $2p-1$ is reduced, Lemma~\ref{lem:elementary}(1) implies that $S$ is not reduced.
Together with $h_p(2)\ge2p-2$, this proves $h_p(2)=2p-2$.
\end{proof}

\begin{remark}\label{rem:peng}
The upper bound $h_p(2)\le2p-1$ also follows from Lemma~\ref{lem:line-bound} and Peng's result $\mathsf c_0(G)=2p-1$.
Assume to the contrary that $S$ is a reduced sequence over $G$ with $|S|\ge2p$.
For any $a\mid S$, the sequence $Sa^{-1}$ is reduced by Lemma~\ref{lem:elementary}(1) and regular by Lemma~\ref{lem:line-bound}.
Since $|Sa^{-1}|=|S|-1\ge2p-1=\mathsf c_0(G)$, the definition of $\mathsf c_0(G)$ implies $\Sigma(Sa^{-1})=G$.
Hence $\Sigma_0(Sa^{-1})=\Sigma_0(S)=G$, whereas~\eqref{eq:reduced} requires $\Sigma_0(Sa^{-1})\subsetneq\Sigma_0(S)$, a contradiction.
Together with the lower bound in the proof of Theorem~\ref{thm:main}, this yields $2p-2\le h_p(2)\le2p-1$.
\end{remark}

\begin{corollary}\label{cor:preserve}
Let $p$ be a prime, let $S$ be a finite sequence over $C_p^2$, and let $Z\mid S$ with $0\nmid Z$ and $|Z|\le p-1$.
There exists a subsequence $R$ such that
\[
 Z\mid R\mid S,\qquad |R|\le2p-2,\qquad \Sigma_0(R)=\Sigma_0(S).
\]
The bound on $|R|$ is sharp even for $Z=1$, and the bound $|Z|\le p-1$ is sharp for every prime $p$.
\end{corollary}

\begin{proof}
Since $S$ has finitely many subsequences and satisfies $Z\mid S$, there exists a subsequence $R\mid S$ of minimum length subject to $Z\mid R$ and $\Sigma_0(R)=\Sigma_0(S)$.
If $0\mid R$, then $0\nmid Z$ implies $Z\mid R0^{-1}$, and $\Sigma_0(R0^{-1})=\Sigma_0(R)$, contrary to the minimality of $|R|$.
Thus $0\nmid R$.

Assume to the contrary that $|R|\ge2p-1$.
Then there exists $T=g_1\cdots g_{2p-1}\mid R$.
We show that at least $p$ indices $i$ satisfy $\Sigma_0(Tg_i^{-1})=\Sigma_0(T)$.
If $T$ is regular, this follows from Theorem~\ref{thm:deletion}, since $\Sigma_0(Tg_i^{-1})=G$ implies $\Sigma_0(T)=G$.
If $T$ is not regular, then $0\nmid T$ implies $|T_H|\ge p$ for some subgroup $H$ of order $p$.
For every $a\mid T_H$, the sequence $T_Ha^{-1}$ contains at least $p-1$ nonzero elements of $H$, so $\Sigma_0(T_Ha^{-1})=H$.
Indeed, if $\varnothing\ne B\subsetneq H$ and $h\in H\setminus\{0\}$, then $B+h\ne B$ since $\langle h\rangle=H$.
Thus $|B\cup(B+h)|\ge |B|+1$ whenever $\varnothing\ne B\subsetneq H$.
Since $\Sigma_0(1)=\{0\}$ and $|H|=p$, successive use of~\eqref{eq:union} shows that the subsequence sums of any $p-1$ nonzero elements of $H$ cover $H$.
By~\eqref{eq:union}, $\Sigma_0(Ta^{-1})=H+\Sigma_0(TT_H^{-1})=\Sigma_0(T)$ for every $a\mid T_H$.
This accounts for at least $|T_H|\ge p$ indices.

Let $V$ be the subsequence of $T$ consisting of the $g_i$ for which $\Sigma_0(Tg_i^{-1})=\Sigma_0(T)$.
Since $|V|\ge p>|Z|$, there exists $a\in G$ with $\mathsf v_a(V)>\mathsf v_a(Z)$.
Then $V\mid R$ implies $a\mid RZ^{-1}$, so $Z\mid Ra^{-1}$.
By~\eqref{eq:union},
\[
 \Sigma_0(Ra^{-1})=\Sigma_0(Ta^{-1})+\Sigma_0(RT^{-1})=\Sigma_0(R)=\Sigma_0(S).
\]
Since $|Ra^{-1}|=|R|-1$, this contradicts the minimality of $|R|$.
Hence $|R|\le2p-2$.

For $Z=1$, the reduced sequence $S_0=e_1^{\,p-1}e_2^{\,p-1}$ in the proof of Theorem~\ref{thm:main} shows that $|R|\le2p-2$ is sharp.
Indeed, for any proper subsequence $R\mid S_0$, there exists $a\mid S_0R^{-1}$, so $\Sigma_0(R)\subseteq\Sigma_0(S_0a^{-1})\subsetneq\Sigma_0(S_0)$.

To see that $|Z|\le p-1$ is sharp for odd $p$, let $S=W$ and let $Z$ be as in~\eqref{eq:deletion-example}, so $|Z|=p$.
If $Z\mid R\mid W$ and $R\ne W$, then $R\mid We_1^{-1}$, since $WZ^{-1}=e_1^{\,p-1}$.
Thus $-e_1\notin\Sigma_0(R)$, although $\Sigma_0(W)=G$.
For $p=2$, let $e_1,e_2$ be a basis of $G$ and let $S=e_1^2e_2$ and $Z=e_1^2$.
If $Z\mid R\mid S$ and $|R|\le2$, then $R=Z$ and $\Sigma_0(R)=\langle e_1\rangle\ne G=\Sigma_0(S)$.
\end{proof}

\section{A group algebra proof of the upper bound}\label{sec:remarks}

Let $p$ be a prime and let $G=C_p^2$.
We prove $h_p(2)\le 2p-2$ using the group algebra $\mathbb F_p[G]$.
For a sequence $U$ over $G$, let $\Sigma(U)$ be the set of nonempty subsequence sums and let $\Sigma_0(U)=\Sigma(U)\cup\{0\}$.
For a subgroup $H\le G$, let $U_H$ denote the subsequence of all elements of $U$ in $H$.
Let $1$ denote the empty sequence, and for $V\mid U$ let $UV^{-1}$ denote the sequence determined by $U=V(UV^{-1})$.
Recall that $U$ is reduced if $\Sigma_0(Ua^{-1})\subsetneq\Sigma_0(U)$ for every $a\mid U$.
The maximum length of such a sequence over $G$ is $h_p(2)$.

The group algebra $\mathbb F_p[G]$ has basis $\{X^g:g\in G\}$, with $X^0=1$ and $X^gX^h=X^{g+h}$.
For $U=u_1\cdots u_\ell$, let
\[
 \Pi(U)=\prod_{i=1}^{\ell}(1-X^{u_i}),\qquad \Pi(1)=1.
\]
Lemma~\ref{lem:group-algebra}(1) follows from Gao, Peng and Zhong~\cite[Lemmas~3.1 and~3.10(1)--(3)]{GPZ}, since the Davenport constant of $C_p^2$ is $\mathsf D(C_p^2)=2p-1$.
Lemma~\ref{lem:group-algebra}(2) is Peng's nonvanishing theorem~\cite[Theorem~1]{PengI} in the form stated in~\cite[Lemma~3.11]{GPZ}.

\begin{lemma}\label{lem:group-algebra}
Let $p$ be a prime and let $U$ be a sequence of nonzero elements of $G=C_p^2$.
\begin{enumerate}
\item If $\Pi(U)\ne0$, then $|U|\le2p-2$.
If, in addition, either $|U|=2p-2$ or $X^g\Pi(U)=\Pi(U)$ for every $g\in G$, then $G\setminus\{0\}\subseteq\Sigma(U)$.
\item Suppose that $p\le |U|\le2p-2$ and $|U_H|\le p-1$ for every subgroup $H$ of order $p$.
If at most $2p-1-|U|$ subgroups of order $p$ satisfy $U_H\ne1$, then $\Pi(U)\ne0$.
\end{enumerate}
\end{lemma}

\begin{proof}[Proof of the upper bound in Theorem~\ref{thm:main}]
Assume to the contrary that $W$ is a reduced sequence over $G$ with $|W|\ge2p-1$, and let $S\mid W$ have length $2p-1$.
Every subsequence $U\mid W$ is reduced.
Indeed, if $\Sigma_0(Ua^{-1})=\Sigma_0(U)$ for some $a\mid U$, then $\Sigma_0(Wa^{-1})=\Sigma_0(Ua^{-1})+\Sigma_0(WU^{-1})=\Sigma_0(W)$, contrary to the definition of a reduced sequence.
In particular, $0\nmid S$, since $\Sigma_0(0)=\Sigma_0(1)$.
For each subgroup $H$ of order $p$, let $S_H=b_1\cdots b_m$.
For $1\le j\le m$, the subsequence $b_1\cdots b_j$ is reduced, so $\Sigma_0(b_1\cdots b_{j-1})\subsetneq\Sigma_0(b_1\cdots b_j)$.
Since $\Sigma_0(1)=\{0\}$ and $\Sigma_0(S_H)\subseteq H$, we have $|S_H|+1=m+1\le|\Sigma_0(S_H)|\le p$.

We construct a subsequence $T\mid S$ with $\Pi(T)\ne0$ by the method used in the proof of~\cite[Theorem~3.12]{GPZ}.
Let $H_1,\ldots,H_t$ be the subgroups of order $p$ such that $S_{H_i}\ne1$, and let $a_i\mid S_{H_i}$ for $1\le i\le t$.
Since $0\nmid S$ and distinct subgroups of order $p$ intersect only in $\{0\}$, we have $a_1\cdots a_t\mid S$.
Let $T=S(a_1\cdots a_t)^{-1}$.
If $|T|\ge p$, then $p\le |T|=2p-1-t\le2p-2$, and $|T_H|\le|S_H|\le p-1$ for every subgroup $H$ of order $p$.
At most $t=2p-1-|T|$ such subgroups satisfy $T_H\ne1$, so $\Pi(T)\ne0$ by Lemma~\ref{lem:group-algebra}(2).

Suppose that $|T|<p$.
Each $u\mid T$ belongs to exactly one subgroup of order $p$, namely $\langle u\rangle$.
Since $G$ has $(p^2-1)/(p-1)=p+1>|T|$ subgroups of order $p$, there exists a subgroup $H$ of order $p$ such that $T_H=1$.
Let $\bar e$ generate $G/H$ and let $z=1-X^{\bar e}$ in $\mathbb F_p[G/H]$.
Since $z^p=1-X^{p\bar e}=0$ and $X^{\bar e}=1-z$ generates $\mathbb F_p[G/H]$, the homomorphism $\mathbb F_p[Z]/(Z^p)\to\mathbb F_p[G/H]$ defined by $Z\mapsto z$ is surjective.
Both algebras have dimension $p$, so $Z\mapsto z$ induces an isomorphism.
Let $\varphi:\mathbb F_p[G]\to\mathbb F_p[G/H]$ be the homomorphism defined by $\varphi(X^g)=X^{g+H}$ for $g\in G$.
For every $u\mid T$, we have $u\notin H$, so $u+H=j\bar e$ for some $1\le j<p$ and $\varphi(1-X^u)=1-(1-z)^j$.
The polynomial $(1-(1-Z)^j)/Z$ has constant coefficient $j\ne0$ in $\mathbb F_p$, so it is relatively prime to $Z^p$ and its residue class modulo $Z^p$ is a unit.
Under the isomorphism $Z\mapsto z$, it follows that $1-(1-z)^j$ is $z$ times a unit in $\mathbb F_p[G/H]$.
Hence $\varphi(\Pi(T))=z^{|T|}v$ for some unit $v\in\mathbb F_p[G/H]$.
Since $|T|<p$, we have $z^{|T|}\ne0$, and therefore $\Pi(T)\ne0$.

Since $\Pi(T)\ne0$, there exists a subsequence $R\mid S$ of maximum length subject to $T\mid R$ and $\Pi(R)\ne0$.
By Lemma~\ref{lem:group-algebra}(1), $|R|\le2p-2$, and $G\setminus\{0\}\subseteq\Sigma(R)$ if $|R|=2p-2$.
If $|R|\le2p-3$, then $SR^{-1}\mid a_1\cdots a_t$ and $|SR^{-1}|\ge2$.
Thus there exist $i\ne j$ such that $a_i a_j\mid SR^{-1}$, and the maximality of $|R|$ implies
\[
 (1-X^{a_i})\Pi(R)=\Pi(Ra_i)=0,\qquad
 (1-X^{a_j})\Pi(R)=\Pi(Ra_j)=0.
\]
Since $\langle a_i\rangle=H_i\ne H_j=\langle a_j\rangle$, the elements $a_i,a_j$ generate $G$.
The equalities $X^{a_i}\Pi(R)=X^{a_j}\Pi(R)=\Pi(R)$ therefore imply $X^g\Pi(R)=\Pi(R)$ for every $g\in G$.
It follows from Lemma~\ref{lem:group-algebra}(1) that $G\setminus\{0\}\subseteq\Sigma(R)$.

In both cases $\Sigma_0(R)=\Sigma(R)\cup\{0\}=G$.
Since $|R|<|S|$, there exists $a\mid SR^{-1}$, so $R\mid Sa^{-1}$ and
\[
 G=\Sigma_0(R)\subseteq\Sigma_0(Sa^{-1})\subseteq\Sigma_0(S)\subseteq G.
\]
Thus $\Sigma_0(Sa^{-1})=\Sigma_0(S)$ although $S$ is reduced, a contradiction.
It follows that $h_p(2)\le2p-2$.
\end{proof}

\section*{Acknowledgements}
This work was supported by the National Natural Science Foundation of China
(Grant No. 12371335 and 12671389).
\section*{Declaration of generative AI and AI-assisted technologies}
The authors used a large language model during the preparation of this manuscript
for assistance with mathematical exploration, checking arguments,
organization and language editing. The authors are responsible for
the mathematical statements, proofs, references and final text.

\begingroup
\makeatletter
\renewcommand{\@biblabel}[1]{#1.}
\makeatother

\endgroup
\end{document}